\documentclass{iise}

\usepackage{multirow}

\conference{Proceedings of the IISE Annual Conference \& Expo 2025 \\
E. Gentry, F. Ju, X. Liu, eds.}% Do not change this line.

\title{\titlesize Dynamic Less-Than-Truckload Transportation Planning in Hyperconnected Hub Networks with Multi-Carrier Operations}

\author{
Tiankuo Zhang, Jingze Li, Benoit Montreuil\\Physical Internet Center, Supply Chain \& Logistics Institute, 
School of Industrial \& Systems Engineering, Georgia Institute of Technology
\\Atlanta, GA}
\authorlist{Zhang, Li and Montreuil}%Heading
\abstractID{5143}% Fill in the web conference management system-assigned abstract ID.

\begin{document}
\maketitle

\begin{abstract}
{\small Less-than-truckload (LTL) shipment is vital in modern freight transportation yet is in dire need of more efficient usage of resources, higher service responsiveness and velocity, lower overall shipping cost across all parties, and better quality of life for the drivers. The industry is currently highly fragmented, with numerous small to medium-sized LTL carriers typically operating within dedicated regions or corridors, mostly disconnected from each other.}

{\small This paper investigates the large-scale interconnection of LTL carriers enabling each to leverage multi-carrier networks for cross-region services exploiting their mutual logistic hubs, in line with Physical Internet principles. In such a network, efficient open cooperation strategies are critical for optimizing multiparty relay shipment consolidation and delivery, transport and logistic operations and orchestration, and enabling inter-hub driver short hauls. }

{\small To dynamically plan relay truck transportation of involved carriers across hyperconnected hub networks, we develop an optimization-based model to build loads, coordinate shipments, and synchronize driver deliveries. We report a simulation-based experiment in a multiparty LTL network covering the eastern U.S. in three scenarios: 1) each carrier operates separately and serves its clients with end-to-end transportation, 2) each carrier operates separately and adopts relay transportation in its service region, and 3) all carriers operate jointly and serve clients in the multi-carrier hyperconnected relay network. By comparing these three scenarios, we evaluate the impact of relay transportation and carrier cooperations on cost savings, trip duration, and greenhouse gas emissions. Overall, this research advances operational efficiencies through an effective collaborative solution across the LTL industry and contributes to the pursuit of sustainable logistics networks.}
\end{abstract}

\section*{Keywords}
Less-than-Truckload Shipping; Multi-Carrier Collaboration; Sustainable Logistics; Relay Truck Transportation; Decision Support Systems

\section{Introduction}
Less-than-truckload (LTL) shipments, in which items from multiple parties (shippers, customers, or consignees) fill a truck or trailer, are the backbone of modern supply chains due to their cost efficiency and fast responsiveness. As more parties are involved, LTL shipments require multiple stops and seamless collaborations, which are less complicated for a single large-scale LTL provider such as UPS and FedEx. With carefully established systems and extensive coverages, large providers can determine efficient plans to navigate each shipment and truck with fast transits. However, these are critical for collaborations between small to medium-sized firms with limited terminals that often operate within dedicated regions. With no collaborations and fragmented service regions, carriers will either entirely deny out-of-region delivery requests or impose extra resources and costs to fulfill such requests. Through collaborations, carriers can jointly handle cross-region shipments with shared resources, including trucks, drivers, and terminals. Nevertheless, complications in consolidation, transshipment, and coordination increase drastically as more parties are involved, potentially leading to high operational costs and significant pressure for on-time services.

Besides the dilemma in carrier collaboration, driver shortages, ranked fourth in the list of supply chain issues by the American Transportation Research Institute (ATRI) in 2023, also contribute to the challenges faced by the LTL industries. According to the American Trucking Association's (ATA) estimation, the U.S. was short by 81258 drivers in 2021. In 2023, with a slightly recovered population, ATA estimated a new shortage of 60000, which remains alarming. One main reason for driver shortage is long working hours, especially in long hauls, which require multi-day trips and allow fewer return-to-domicile opportunities, resulting in drivers living in poor conditions with health concerns. Traditional supply chain operations also lead to low truck utilization, with trailers approximately 60\% complete when traveling loaded in the U.S. and less than 10\% global efficiency. Trucks also become less loaded as they travel and drop off goods, causing even lower overall utilization rates \cite{b01}. Sustainability is another concern of modern supply chains. As transport accounts for 25\% of global carbon dioxide emissions, many countries and groups have set goals in pursuing net-zero supply chains \cite{gmr01}.

In response to the challenges above, this study proposes a decision support tool to effectively plan collaborative LTL shipments for regional small to medium-sized LTL companies and provide drivers with more return-to-domicile opportunities. The decision support tool considers a collaborative hyperconnected relay LTL hub network inspired by Physical Internet (PI) principles, enabling cross-region collaboration and resource sharing among carriers. To orchestrate such a system, we developed an optimization-based model to build loads, coordinate shipments, and synchronize driver deliveries. We incorporated relay transportation to allow drivers to return to their base terminal while ensuring shipment timeliness. In three scenarios, we apply our system to a simulation-based experiment of a multi-carrier collaboration in eastern U.S. states. The results iteratively validate PI principles' application in the real world and show our study's effectiveness in improving cost-efficiencies, allowing drivers to return to their domiciles, and reducing greenhouse gas emissions. 

The rest of the paper is structured as follows. Section 2 presents the related literature. Section 3 proposes our model formulation and potential variants. Section 4 analyzes our study's effectiveness through a simulation-based experiment of two scenarios. Section 5 summarizes the contributions, limitations, and future work directions.

\section{Literature Review}
Researchers have studied carrier collaborations and their effects through various empirical analyses on productivity, market positioning, and information sharing. Many studies focused on the overall increment in profitability and impacts on individual gains, and Cruijssen et al. claimed a net gained value of 30 percent \cite{kk01, cdf01, cbdfs01}. More recently, Ferrell et al. identified the lack of effective collaborations in practice and explored open collaboration opportunities between logistics providers with PI principles \cite{fekr01}.

In mathematical modeling, many previous studies have investigated collaborative truck transportation as collaborative pickup and delivery problems with centralized solution models based on mixed integer or heuristic models \cite{GANSTERER20181}. In these settings, there are no restrictions on service region, and all carriers have equal access to all pickup points. Other researchers studied the application of service network design in collaborative truck transportation as collaborative service network design (CSND) problems. Initially defined by studies in liner shipping, each carrier has a service network region in these service networks with time-sensitive and cross-region delivery requests \cite{ae01}. Zhang et al. studied the CSND problem, focusing on disjoint service regions and demand uncertainty \cite{zhw01}. Yet, previous research on CSND problems considered neither limited vehicle nor vehicle balance constraints, and our study fills this research gap.

Researchers have conducted comprehensive research on concepts of PI, especially the topic of hyperconnected transportation planning and its economic, environmental, and social efficiency impacts \cite{hmsbp01}. Through modeling and planning, Orentein and Raviv applied hyperconnected networks to small parcel delivery, involving intermediate service points and exploiting consolidation opportunities \cite{or01}. Li et al. (2022) introduced relay transportation and developed an operating system in hyperconnected networks with multi-agent architecture to generate and coordinate shipments, plan tractors and trailers, and schedule truckers, with the potential to increase drivers' return-home opportunities \cite{li01}. Li et al. then studied the combination of stochastic models and hyperconnected relay hub networks with a two-stage solution model based on service network design \cite{li02}. Previous research on PI-related hyperconnected networks focused on improvements to a single logistics provider or the overall system, and benefits to each participating carrier are less emphasized. This research fills the gap by validating improvements on each carrier's logistics system by applying PI principles in a collaborative setting.

\section{Methodology}
\begin{figure}[htb]
    \centering
    \includegraphics[width=6in]{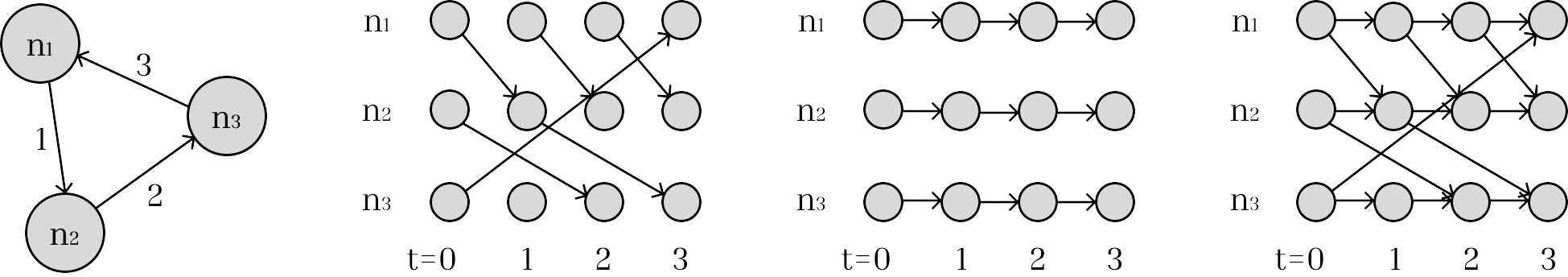}
    \caption{Example of $\mathcal{G}^P$, $\mathcal{A}$, $\mathcal{H}$, and $\mathcal{G}$ with $T=3$}\label{fig1}
\end{figure}

Our study considers a centralized platform that plans logistics decisions for each carrier based on a given multi-carrier hub network. The network is denoted by $\mathcal{G}^P = (\mathcal{N}^P, \mathcal{A}^P)$, where $\mathcal{N}^P$ denotes the hub node set, and $\mathcal{A}^P$ denotes the directed arc set connecting hub node pairs with respective transit time $\tau$. $\mathcal{T} = {0, 1, ..., T}$ represents the planning instances starting from instance 0 and ending at instance T. Based on $\mathcal{T}$, we discretize the network $\mathcal{G}^P$ and obtain the time-expanded network $\mathcal{G} = (\mathcal{N}, \mathcal{A}\cup\mathcal{H})$. In the time-expanded network, $\mathcal{N} = \{(n, t)|n\in\mathcal{N}^P,t\in\mathcal{T}\}$ are copies of $\mathcal{N}^P$ at each planning instance; $\mathcal{A} = \{((n_1, t), (n_2, t + \tau))|(n_1, t), (n_2, t + \tau)\in\mathcal{N}\}$ are moving arcs are copies of $\mathcal{A}^P$ connecting node copies considering their transit time $\tau$; $\mathcal{H}$ are holding timed arcs connecting each node copy to the next consecutive copy, i.e., $(n_1, t)$ and $(n_1, t + 1)$. For simplicity, we refer to node copies in $\mathcal{N}$ as timed nodes and arc copies and holding arcs in $\mathcal{A}\cup\mathcal{H}$ as timed arcs. Figure \ref{fig1} demonstrates an example of $\mathcal{G}^P$ with three nodes and three arcs and the time-expanded $\mathcal{G}$, when $T = 3$. 

\hspace{0pt}$\mathcal{G}^P$\hspace{100pt}$\mathcal{A}$\hspace{110pt}$\mathcal{H}$\hspace{110pt}$\mathcal{G}$

The centralized platform plans truck routings based on given commodity orders $k\in\mathcal{K}$ with given origin $o_k$, destination $d_k$, release time $r_k$, deadline $l_k$, and volume $v_k$. $c^k_a$ denotes the variable cost of transporting in on timed arc $a\in\mathcal{A}$, which captures the difference in commodity type, transshipment methods, and experience level of the driver. The multi-carrier model considers a set of carriers $\mathcal{S}$, each with a given service region encompassing a set of hub nodes that are not necessarily disjoint and denoted by $\{\mathcal{N}_s\subseteq\mathcal{N}|s\in\mathcal{S}\}$. The platform also considers vehicle routing such that the number of trucks is limited with vehicle conservation. For a carrier $s\in\mathcal{S}$ $m_{sn}^0$ denotes the given number of s' trucks starting at hub node $n\in\mathcal{N}^P$, and $m_{sn}^T$ denotes the minimum number of s' trucks required to end at hub node $n\in\mathcal{N}^P$. We consider uniform trucks with maximum volumes $v^{max}$ and fixed cost $f_a$ for each truck transporting commodities on $a\in\mathcal{A}$. 

The mathematical formulation with two variables is shown below. $X_{ka}$ is a binary variable indicating whether a commodity $k\in\mathcal{K}$ travels on timed arc $a \in \mathcal{A}$. $Y_{sa}$ is an integer variable representing the number of trucks of carrier $s\in S$ allocated on timed arc $a \in \mathcal{A}$. $\mathcal{A}^+(n,t)$ and $\mathcal{H}^+(n,t)$ denote the set of timed moving and holding arcs leaving a timed node $(n, t)$. $\mathcal{A}^-(n,t)$ and $\mathcal{H}^-(n,t)$ denote the set of timed moving and holding arcs entering a timed node $(n, t)$. The objective function (1) minimizes the sum of the total variable cost of each commodity and the total fixed cost for allocating trucks on each timed arc. Constraint (2) ensures proper transfer of commodities such that each commodity is released at its release time by the origin, is due at its deadline by the destination, and flows through other timed nodes with conservation. Constraint (3) allocates trucks on each timed arc based on commodity shipments. Constraint (4) sets the initial number of trucks of each carrier at each hub, while Constraint (5) sets the final lower bounds. Constraint (6) enforces vehicle conservation.

\begin{align}
min&\sum_{k \in \mathcal{K}}\sum_{a \in \mathcal{A}}c_a^kX_{ka} + \sum_{s\in S}\sum_{a \in \mathcal{A}}f_aY_{sa}\\
    s.t.&\sum_{a\in\mathcal{A}^+(i) \cup \mathcal{H}^+(i)} X_{ka} - \sum_{a\in\mathcal{A}^-(i)\cup \mathcal{H}^-(i)} X_{ka} = \begin{cases}
    1, & \text{if } i = (o_k, r_k)\\
    -1, & \text{if } i = (d_k, l_k)\\
    0, & \text{otherwise}\\
    \end{cases}\hspace*{20pt}\hfill,\forall k \in \mathcal{K}\\
    &\sum_{k\in\mathcal{K}} X_{ka}\leq \sum_{s\in\mathcal{S}}v^{max}Y_{sa}, \forall a\in\mathcal{A}\\
    &\sum_{a\in\mathcal{A}^+(n, 0)\cup \mathcal{H}^+(n, 0)}Y_{sa} = m_{sn}^0, \forall s\in\mathcal{S}, n \in \mathcal{N}\\
    &\sum_{a\in\mathcal{A}^-(n, T)\cup \mathcal{H}^+(n, T)}Y_{sa} \geq m_{sn}^T, \forall s\in\mathcal{S}, n \in \mathcal{N}\\
    &\sum_{a\in\mathcal{A}^+(n, t)\cup \mathcal{H}^+(n, t)}Y_{sa} = \sum_{a\in\mathcal{A}^-(n, t)\cup \mathcal{H}^-(n, t)}Y_{sa}, \forall s\in\mathcal{S}, n \in \mathcal{N}, t\notin\{0, T\}
\end{align}

\section{Results and Discussion}
\begin{figure}[htb]
    \centering
    \includegraphics[width=5.4in]{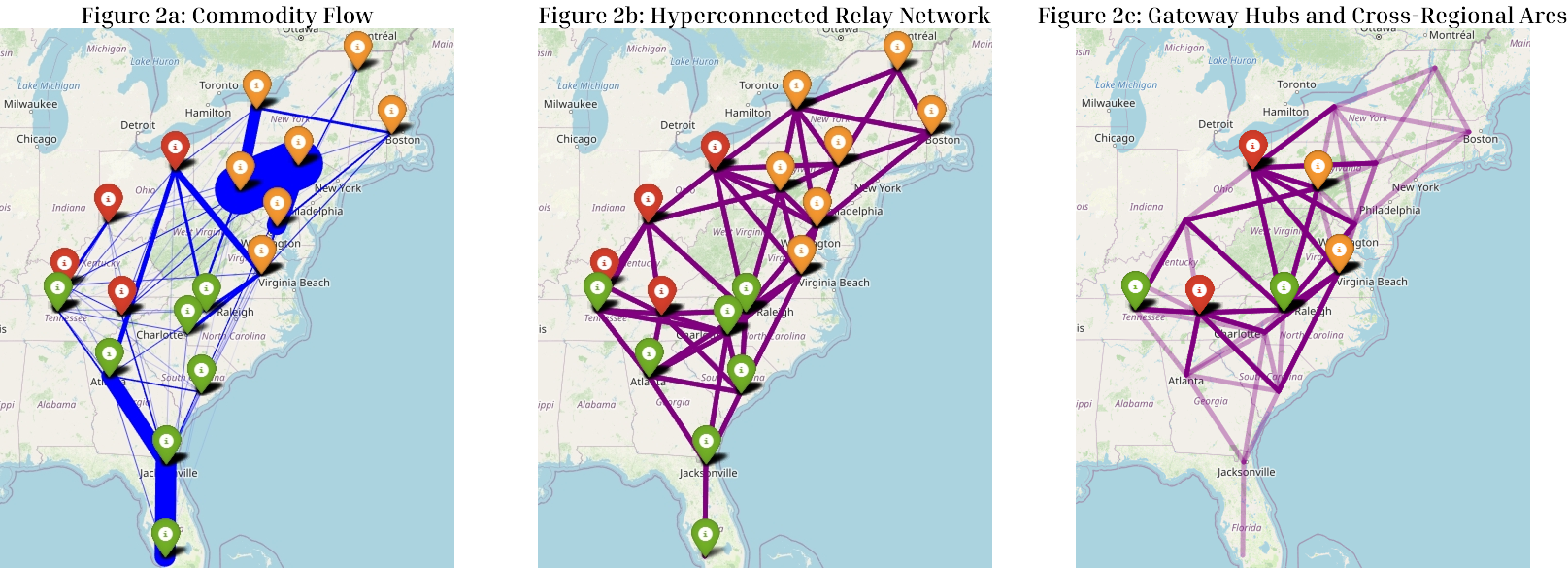}
    \caption{Commodity Flow (Blue) and Hyperconnected Network (Purple)}\label{fig2}
\end{figure}

We present a simulation-based experiment pooled from a hub network of a real-life national alliance between four LTL carriers, and we focus on their operations involving three carriers in eastern U.S. states. The representative hub network includes 18 hubs and 120 arcs ($< 5.5$ hours). We consider a time horizon of 48 hours with 6 extra winding-down hours for returning vehicles. We configure hourly planning instances and simulate 82 commodity requests based on the Freight Analysis Framework (FAF) data of the Bureau of Transportation Statistics. We assume each request is consolidated from more granular shipper-specific requests with the exact origin, destination, and a similar time window.  Figure \ref{fig2}a shows the derived commodity flow in blue, with the colors of the hubs indicating the belonging carriers. Notably, there are many cross-regional requests.

We compare our model results between three scenarios: (1) all carriers operate independently with only end-to-end transportation, (2) all carriers operate independently and adopt relay transportation in their respective service region, and (3) all carriers operate in the collaborative hyperconnected relay network. In scenario 1, each carrier independently solves the introduced model with end-to-end routing; thus, all commodities will travel directly from origin to destination with consolidation opportunity only at their origins. In scenario 2, each carrier solves a similar model separately with in-region relay transportation while all in-region arcs are less than 5.5 hours. This ensures that every truck and its driver will complete a round trip within 11 hours, abiding by the federal hour-of-service regulation. Under this scenario, carriers can consolidate each commodity at multiple stops within their service region. In scenario 3, we propose the collaborative hyperconnected relay network, and the centralized platform solves the joint model with relay transportation. This allows each commodity to be consolidated at all hubs. Figure \ref{fig2}b demonstrates the hyperconnected LTL hub network employed in this scenario. In addition, we only allow cross-regional arcs leaving from a preselected group of gateway hubs, highlighted in Figure \ref{fig2}c. 

\begin{figure}[htb]
    \centering
    \includegraphics[width=6.5in]{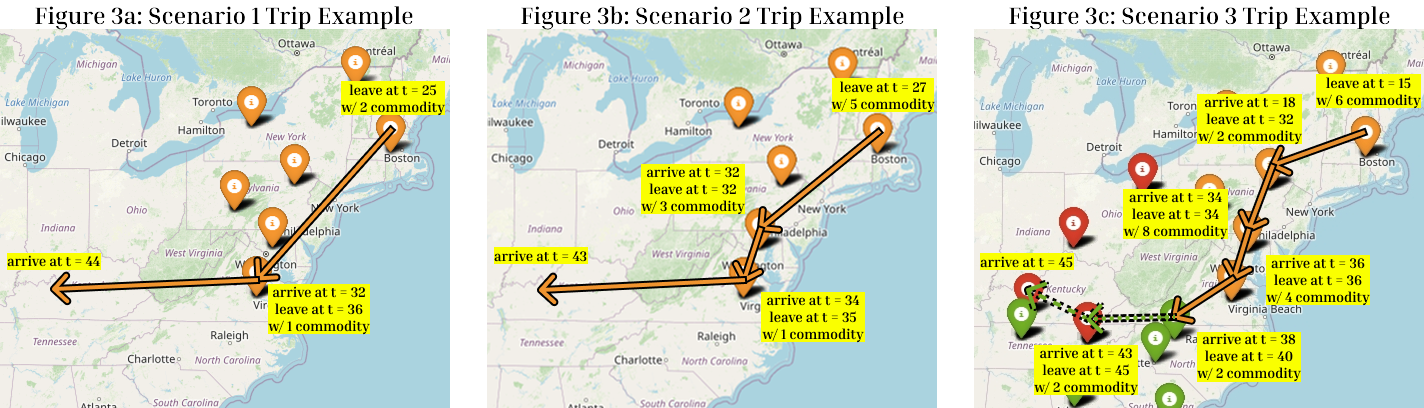}
    \caption{Optimized Truck Trip Example of All Scenarios}\label{fig3}
\end{figure}

Based on our tool's optimized results in each scenario, Figure \ref{fig3} shows an exemplary truck trip, which starts and ends at the same hub nodes and fulfills a cross-regional commodity request. In the first scenario, the truck leaves its starting hub and delivers two commodities to its first stop, from which it delivers a cross-regional commodity. In the second scenario, due to the implementation of an in-region relay network, the truck conducts more trips in its carrier's service region with more commodities and, eventually, delivers the cross-regional commodity. In the first two scenarios, trucks from the orange carrier have no access to other carriers' hubs, which are thus not shown on the map. In the third scenario, with the collaborative hyperconnected relay network, the truck takes trips with even more commodities. It leaves its gateway hub and delivers its cross-regional commodities to an out-of-region hub, from which another carrier's trucks (indicated by the dashed arrows and different colors) will take over the delivery task of the cross-regional commodity. These differences validate our tool's effectiveness in increasing truck utilization rates and improving drivers' working conditions. Drivers in the third scenario do not need to leave their service region for long trips in the first and second scenarios, while commodities will still be delivered on time.

We also compare the three scenarios on five key performance indicators (KPIs): cost, total truck hours, total truck trips, average trip hour, and emission. In assessing cost, we utilized data from the 2023 trucking operational cost analysis of ATRI and the 2024 annual freight trucking rate report of Uber Freight. We assumed that regular truck drivers handle in-region truck trips, and cross-region truck trips require dedicated truck drivers with higher per-mile wages. Trips over 5.5 hrs also induce extra costs, accommodating drivers' on-the-road expenses, such as food and housing. In estimating greenhouse gas emissions, we adopt Environmental Defense Funds' ton-mile estimation function, encompassing both distance and weight in greenhouse gas emission calculation. In the collaborative scenario, carriers will share costs on shared trucks based on weight ratios when allocating costs and emissions.

\begin{figure}[htb]
    \centering
    \includegraphics[width=5.1in]{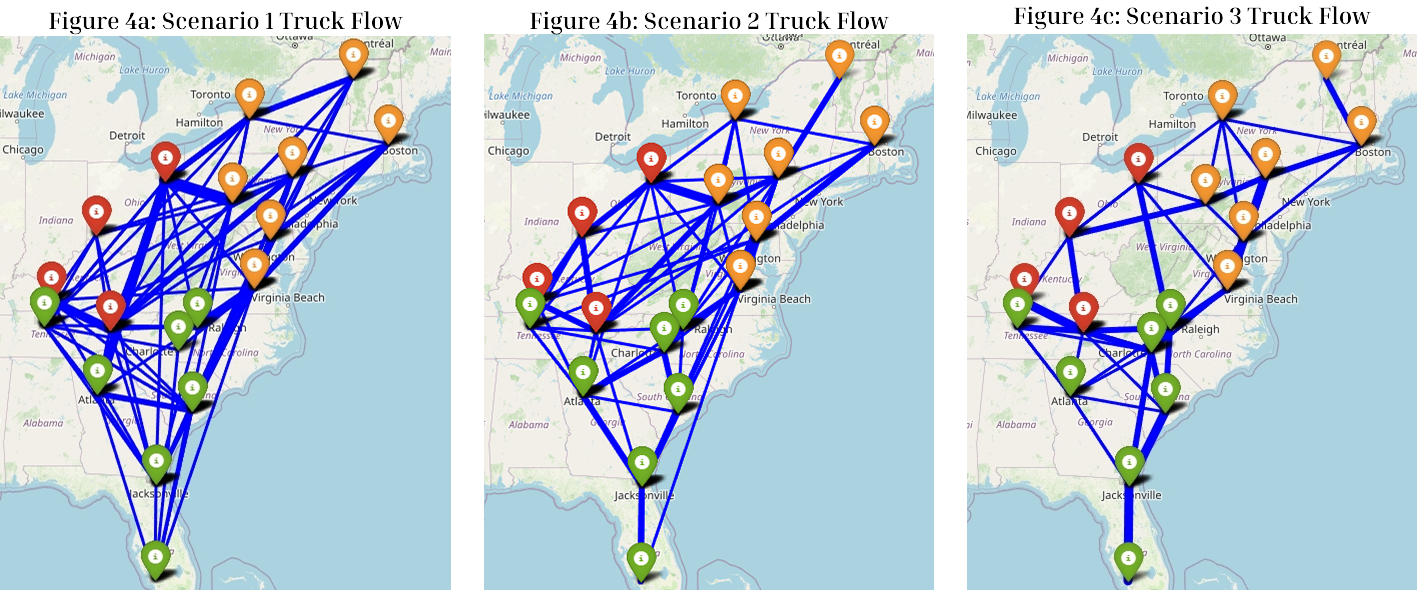}
    \caption{Truck Flow Comparison Between All Scenarios}\label{fig4}
\end{figure}

Table \ref{tab1} summarizes our results by carriers in each scenario. Comparing the first and the second scenarios, applying in-region relay transportation reduced costs by 26\%, total trip hours by 30\%, and emissions by 5\%. This difference validates relay transportation's effectiveness in cost-efficiency, drivers' satisfaction, and sustainability. Yet, the second scenario requires more trucking trips since some in-region trips exceeding 5.5 hours are no longer feasible. The collaborative hyperconnected relay network implemented in the third scenario excelled in most KPIs of all carriers. It reduced the overall cost (53\% of the first scenario and 72\% of the second scenario) and transportation hours (44\% of the first scenario and 63\% of the second scenario). The total number of trips decreases less drastically, and carrier A needs to conduct one more trip in the third scenario than the first scenario. Due to its construction, the collaborative hyperconnected relay network enables all truck drivers to return to their domiciles daily with lower average trip hours. This is supported by the truck flow comparison shown in Figure \ref{fig4}, which shows a gradual decrease in trip density. Carriers will also produce total greenhouse emissions of 89\% of the first scenario. Our proposed network's significant improvements demonstrate the importance of collaboration in logistics, aligning with the PI's core ideology.

\begin{table}[htb]
\caption{Experiment Comparison Results}\label{tab1}
\vspace{-0.7cm}
\begin{center}
\begin{tabular}{l|ccc|ccc|ccc}
\hline
{} & \multicolumn{3}{c|}{\textbf{Scenario 1}} & \multicolumn{3}{c|}{\textbf{Scenario 2}}& \multicolumn{3}{c}{\textbf{Scenario 3}}\\ \hline
\textbf{Carriers} & A & B & C & A & B & C & A & B & C\\ \cline{1-10}
\textbf{Cost (\$)} & 28670 & 57473 & 36834 & 23276 & 42434 & 25176 & 14073 & 31311 & 20277
\\ \cline{1-10}
\textbf{Total Hours} & 131 & 260 & 190 & 105 & 181 & 125 & 87 & 97 & 75
\\ \cline{1-10}
\textbf{Total Trips} & 29 & 42 & 29 & 30 & 38 & 34 & 30 & 31 & 25
\\ \cline{1-10}
\textbf{Average Trip Hour} & 4.52 &  6.19 & 6.55 & 3.5 & 4.76 & 3.68 & 2.9 & 3.13 & 3
\\ \cline{1-10}
\textbf{Emission (metric-ton)} & 95.02 &  203.34 & 132.41 & 91.37 & 192.02 & 122.78 & 85.86 & 182.33 & 116.34
\\ \cline{1-10}
\end{tabular}
\end{center}
\end{table}

\section{Conclusion}
Our study has three main contributions. First, it proposes a decision support tool based on the collaborative hyperconnected relay logistics network to advance cost efficiency, improve drivers' working conditions, and promote sustainability. Second, it provides a mathematical optimization formulation that builds loads, coordinates shipments, and routes vehicles in a centralized platform with application to multiple scenarios. Third, it conducts a simulation-based experiment based on a real-life national alliance between four LTL carriers and cost analyses conducted by ATRI and Uber Freight.

This research also opens several avenues for future work. The first direction is to consider different cost allocation systems to capture various incentive functions of carriers in real life. The second direction is to develop more advanced optimization models to schedule logistics plans in large networks with more hub nodes, transportation arcs, and commodity requests. The third direction involves Well-to-Wheel emission calculations to capture greenhouse gas produced in each trip conclusively.

\end{document}